\documentclass[11pt]{article}
\usepackage{amsmath}
\begin{document}
\begin{center}
\Large{On extensions of Lie algebras}\\
\footnote{Email for correspondence: lsimonian@gmail.com}
\large{L. A. Simonian}\\
\vspace{.4cm}
\large{Abstract}
\end{center}

In the note some construction of Lie algebras is introduced. It is proved that the construction has the same property as a well known wreath product of groups [1]: Any extension of groups can be embedded into their wreath product [2].\\
\bigskip

Let $M$ and $L$ be Lie algebras over an arbitrary field $K$, $U=U\left( L\right)$ - a universal enveloping algebra
of Lie algebra $L$, $\left\lbrace e_{i}, i\in I \right\rbrace$ - a well-ordered basis in $L$. We can
convert the linear space $Hom_{K}\left( U,M\right)$ into a Lie algebra if we define a Lie product by
Leibniz formula $$\left[ f,h\right]  \left( E\right) =\sum_{I*J=E}\left[ f\left( I\right) ,h\left( J\right) \right],$$ where $E, I, J$ are standard monomials [3]:
for example $E=e_{j}e_{k}\cdots e_{m}e_{n}$, $e_{j}\leq e_{k}\cdots \leq e_{m}\leq e_{n}$, $I*J$ is a product
 in a symmetrical algebra of $L$ and $f,h$ are elements of $Hom_{K}\left( U,L\right)$.
 
Define an action of $L$ on  $Hom_{K}\left( U,M\right)$ by a rule $$\left( fu\right) \left( E\right) =
f\left( uE\right),$$ where $f\in Hom_{K}\left( U,M\right)$, $u\in L$ and the product of $u$ and $E$ is taken in algebra $U$.

It can be immediately checked that $Hom_{K}\left( U,M\right)$ is indeed a Lie algebra with respect to the above defined product
 and $L$ acts on $Hom_{K}\left( U,M\right)$ as a Lie algebra of derivations of the Lie algebra $Hom_{K}\left(U,L\right)$.
 
We denote a semidirect product of $Hom_{K}\left( U,M\right)$ and $L$ by $M\;Wr\; L$ and call it a wreath 
product of Lie algebras $M$ and $L$.

The notation and the name are justified by the following Therem that we prove here: 

\textbf{Any extension $ N $ of Lie algebra $ M $ by Lie algebra $ L $ can be embedded into their wreath product $ M\;Wr\;L $}

The theorem is similar to the well known theorem of Kaloujnine and Krasner [2].

To prove it, we need the following way of constructing the extension $N$ via a factor set $g\left( u,v\right)$.

Let $ M $ and $ L $ be Lie algebras and suppose that elements of $L$ act on $M$ as derivations of algebra $M$,
that is $$[x,y]u=[xu,y]+[x,yu].$$ 
Let $g:L\times L\rightarrow M$ be a bilinear mapping, such that
\begin{align}
\mbox{(a)}\; & g(u,v)=-g(v,u) \nonumber \\
\mbox{(b)}\; & g(u,v)w+g([u,v],w)+g(v,w)u+g([v,w],u)+  \nonumber \\
 &    g(w,u)v+g([w,u],v)=0  \nonumber \\
\mbox{(c)}\; & (xu)v-(xv)u=x[u,v]+[x,g(u,v)]  \nonumber
\end{align}
where $x,y\in M$, $u,v\in L$.
Then the direct product $N=M\times L$ of linear spaces $M$ and $L$ can be converted into a Lie algebra by the formula
$$[(x,u),(y,v)]=([x,y]+xv-yu+g(u,v),[u,v]).$$
It can be verified that $N$ is the extension of $M$ by  $L$ with a given factor set $g(u,v)$ and given an action of elements of $L$ on $M$.

We will henceforth assume that $N$, as the extension of $M$ by $L$, is given as just described.

We are now coming to the proof of the Theorem. We will construct an embedding $ \varphi:N\rightarrow M\;Wr\;L $.

If $(x,u)\in N$, then $(x,u)=(x,0)+(0,u)$. Therefore it is enough to determine $\varphi((x,0))$ and
$\varphi((0,u))$. In turn, if $\left\lbrace z_{q}, q\in Q\right\rbrace$ is a basis in $M$ and $$x=\sum_{q}\beta_{q}z_{q},$$ then $\varphi((x,0))$ must equal $\sum_{q} \beta_{q}\varphi((z_{q},0)).$
Therefore it suffices to determine $\varphi((z_{q},0))$. Equally, to determine $ \varphi((0,u)) $ we need to know $\varphi((0,e_{i}))$.

Next, set $\varphi((x,u))=( f_{(x,u)} ,u)$, where $f_{(x,u)}\in Hom_{K} \left( U,M\right)$. In the same 
sense we will use notations $f_{(x,0)},\;f_{(0,u)},\;f_{(z_{q},0)},\; f_{(0,e_{i})}$. For example, $\varphi((x,0))=(f_{(x,0)},0)$ and $\varphi((0,e_{i}))=(f_{(0,e_{i})},e_{i})$.

If $$\varphi([(x,u),(y,v)]=[\varphi((x,u)),\varphi((y,v))],$$ then
\begin{align}
\mbox{(1)}\; & f_{(g(u,v),0)} +f_{(0,[u,v])}=[f_{(0,u)},f_{(0,v)}]+f_{(0,u)}v-f_{(0,v)}u  \nonumber \\
\mbox{(2)}\; & f_{(xu,0)}=[f_{(x,0)},f_{(0,u)}]+f_{(x,0)}u \nonumber \\
\mbox{(3)}\; & f_{([x,y],0)} =[f_{(x,0)},f_{(y,0)}]  \nonumber
\end{align}
and vice versa. 

Now, determine $f_{(z_{q},0)}$ and $f_{(0,e_{i})}$ on standard monomials $E$ by induction in such a way
that assures (1), (2), (3).

Put $f_{(x,0)}(1)=x$ and $f_{(0,u)}(1)=0$.

 If $E=e_{j}$ then for $u=e_{i}$ and $v=e_{j}$, (1) gives us:
$$f_{(g(e_{i},e_{j}),0)}(1) +f_{(0,[e_{i},e_{j}])}(1)=[f_{(0,e_{i})},f_{(0,e_{j})}](1)+(f_{(0,e_{i})}e_{j})(1)-(f_{(0,e_{j})}e_{i})(1) $$
or  $$ g(e_{i},e_{j})=f_{(0,e_{i})}(e_{j})-f_{(0,e_{j})}(e_{i}).$$

Define $f_{(0,e_{i})}(e_{j})$ in the form $\alpha g(e_{i},e_{j})$ where $\alpha\in K$ is to be 
determined.
Then $$f_{(0,e_{j})}(e_{i})=  \alpha g(e_{j},e_{i})= -\alpha g(e_{i},e_{j})=- f_{(0,e_{i})}(e_{j})$$
and $g(e_{i},e_{j})=2f_{(0,e_{i})}(e_{j})$. Hence $$f_{(0,e_{i})}(e_{j})=\frac{1}{2}g(e_{i},e_{j}).$$

To determine $f_{(z_{q},0)}(e_{j})  $ we use (2) and put $ u=e_{j} $ and $ x=z_{q} $:
$$ f_{(z_{q}e_{j},0)}(1)=[f_{(z_{q},0)},f_{(0,e_{j})}](1)+(f_{(z_{q},0)}e_{j})(1). $$ So $ f_{(z_{q},0)}
(e_{j})=z_{q}e_{j}. $
It is immediate that  $$ f_{(x,0)}(e_{j})=xe_{j}. $$

If $ E=e_{j}e_{k} $ then we use
$$ f_{(g(e_{i},e_{j}),0)}(e_{k}) +f_{(0,[e_{i},e_{j}])}(e_{k})=[f_{(0,e_{i})},f_{(0,e_{j})}]
(e_{k})+(f_{(0,e_{i})}e_{j})(e_{k})-(f_{(0,e_{j})}e_{i})(e_{k}).$$ By the previous $ f_{(g(e_{i},e_{j}),0)}(e_{k})=g(e_{i},e_{j})e_{k} $. Next if $$ 
[e_{i},e_{j}]=\sum_{r}\alpha_{r}e_{r} $$ then $$ f_{(0,[e_{i},e_{j}])}(e_{k})=\sum_{r}\alpha_{r}f_{(0,e_{r})}(e_{k}) $$
and the values $ f_{(0,e_{r})}(e_{k}) $ are already known.
We have also $$ [f_{(0,e_{i})},f_{(0,e_{j})}](e_{k})=[f_{(0,e_{i})}(e_{k}),f_{(0,e_{j})}(1)]+[f_{(0,e_{i})}(1),f_{(0,e_{j})}(e_{k})]=0 ,$$
$$ (f_{(0,e_{i})}e_{j})(e_{k})=f_{(0,e_{i})}(e_{j}e_{k}), $$  $$ (f_{(0,e_{j})}e_{i})(e_{k})=f_{(0,e_{j})}(e_{i}e_{k}).$$

If $ e_{i} \leq e_{k}$ then we set $f_{(0,e_{j})}(e_{i}e_{k})=-f_{(0,e_{i})}(e_{j}e_{k}) . $
Then$$ f_{(0,e_{i})}(e_{j}e_{k})=\frac{1}{2}( f_{(g(e_{i},e_{j}),0)}(e_{k}) +f_{(0,[e_{i},e_{j}])}(e_{k})). $$

 In the case of $ e_{i} > e_{k} $ we have
$$f_{(0,e_{j})}(e_{i}e_{k})=f_{(0,e_{j})}(e_{k}e_{i})+f_{(0,e_{j})}([e_{i},e_{k}]) $$ and $ e_{j}\leq e_{k}< e_{i} $.
We set as before $ f_{(0,e_{k})}(e_{j}e_{i})=-f_{(0,e_{j})}(e_{k}e_{i}) $ in
$$ f_{(g(e_{j},e_{k}),0)}(e_{i}) +f_{(0,[e_{j},e_{k}])}(e_{i})=$$ $$[f_{(0,e_{j})},f_{(0,e_{k})}]
(e_{i})+(f_{(0,e_{j})}e_{k})(e_{i})-(f_{(0,e_{k})}e_{j})(e_{i}).$$
Then $$ (f_{(0,e_{j})}(e_{k}e_{i})=\frac{1}{2}( f_{(g(e_{j},e_{k}),0)}(e_{i}) +f_{(0,[e_{j},e_{k}])}(e_{i})) .$$
This imlies
 $$ f_{(0,e_{i})}(e_{j}e_{k})= f_{(g(e_{i},e_{j}),0)}(e_{k}) +$$ $$f_{(0,[e_{i},e_{j}])}(e_{k})+\frac{1}{2}( 
 f_{(g(e_{j},e_{k}),0)}(e_{i}) +f_{(0,[e_{j},e_{k}])}(e_{i})) + f_{(0,e_{j})}([e_{i},e_{k}]).$$
 
Next we determine $f_{(z_{q},0)}(e_{j}e_{k})  $. We set $ u=e_{j} $ and $ x=z_{q} $ in (2). We have
 $$ f_{(z_{q}e_{j},0)}(e_{k})=[f_{(z_{q},0)},f_{(0,e_{j})}](e_{k})+(f_{(z_{q},0)}e_{j})(e_{k}) $$ or
 $$ f_{(z_{q}e_{j},0)}(e_{k})=[f_{(z_{q},0)}(1),f_{(0,e_{j})}(e_{k})]+f_{(z_{q},0)}(e_{j}e_{k}) $$ or
 $$ f_{(z_{q},0)}(e_{j}e_{k})=z_{q}e_{j}e_{k} -\frac{1}{2}[z_{q},g(e_{j},e_{k})] .$$
It follows immediately, that$$ f_{(x,0)}(e_{j}e_{k})=xe_{j}e_{k} -\frac{1}{2}[x,g(e_{j},e_{k})]. $$

Suppose now that $ f_{(x,0)} $ and $ f_{(0,e_{i})} $ are already defined for any standard monomial of degree less 
than $ n $ and let $ E=e_{j}e_{k}F $ be a standard monomial of degree $ n $.

It has to be by (1)
$$ f_{(g(e_{i},e_{j}),0)}(e_{k}F) +f_{(0,[e_{i},e_{j}])}(e_{k}F)=[f_{(0,e_{i})},f_{(0,e_{j})}]
(e_{k}F)+$$ $$(f_{(0,e_{i})}e_{j})(e_{k}F)-(f_{(0,e_{j})}e_{i})(e_{k}F)$$ or
$$ (f_{(0,e_{i})}(e_{j}e_{k}F)-(f_{(0,e_{j})}(e_{i}e_{k}F)=f_{(g(e_{i},e_{j}),0)}(e_{k}F) +f_{(0,
[e_{i},e_{j}])}
(e_{k}F)-$$ $$[f_{(0,e_{i})},f_{(0,e_{j})}](e_{k}F). $$

If $ e_{i} \leq e_{k} $  we put  $f_{(0,e_{j})}(e_{i}e_{k}F)=-f_{(0,e_{i})}(e_{j}e_{k}F)  $ and obtain 
$$ f_{(0,e_{i})}(E)=\frac{1}{2}( f_{(g(e_{i},e_{j}),0)}(e_{k}F) +f_{(0,[e_{i},e_{j}])}
(e_{k}F)-[f_{(0,e_{i})},f_{(0,e_{j})}](e_{k}F )). $$
 
 In the case of $ e_{i} > e_{k} $ we have
 $$ e_{i}e_{k}F=e_{k}G+\sum_{s}\alpha_{s}H_{s}, \; \alpha_{s}\in K. $$ Here $ G $ is a standard monomial which is 
equal to the product of  $e_{i}$ and all factors of $ F $ and $ H_{s} $ are standard monomials of degree less than 
$n$.The first factor of $ G $ can be $ e_{i} $ or the first factor $ e_{m} $ of $ F $. We have $ e_{j}\leq e_{k}< 
e_{i} $ for the first case and $  e_{j}\leq e_{k}\leq e_{m} $ for the second one.

So we set as before $ f_{(0,e_{k})}(e_{j}G)=-f_{(0,e_{j})}(e_{k}G) $ in
$$ f_{(g(e_{j},e_{k}),0)}(G) +f_{(0,[e_{j},e_{k}])}(G)=[f_{(0,e_{j})},f_{(0,e_{k})}]
(G)+f_{(0,e_{j})}(e_{k}G)-f_{(0,e_{k})}(e_{j}G).$$
Then  $$ f_{(0,e_{j})}(e_{k}G)=\frac{1}{2}( f_{(g(e_{j},e_{k}),0)}(G) +f_{(0,
[e_{j},e_{k}])}(G)-[f_{(0,e_{j})},f_{(0,e_{k})}](G)) $$
Ultimately we have
$$ f_{(0,e_{i})}(E)= f_{(g(e_{i},e_{j}),0)}(e_{k}F) +f_{(0,[e_{i},e_{j}])}
(e_{k}F)-[f_{(0,e_{i})},f_{(0,e_{j})}](e_{k}F ) +$$  $$\frac{1}{2}( f_{(g(e_{j},e_{k}),0)}(G) +f_{(0,
[e_{j},e_{k}])}(G)-[f_{(0,e_{j})},f_{(0,e_{k})}](G))+\sum_{s}\alpha_{s}f_{(0,e_{j})}(H_{s}),$$
and values of functions in the right hand side are known.

To determine $f_{(z_{q},0)}(E)  $ we use (2) and put $ u=e_{j} $ and $ x=z_{q} $.
We have
$$ f_{(z_{q}e_{j},0)}(e_{k}F)=[f_{(z_{q},0)},f_{(0,e_{j})}](e_{k}F)+(f_{(z_{q},0)}e_{j})(e_{k}F). $$
This implies
$$f_{(z_{q},0)}(E)= f_{(z_{q}e_{j},0)}(e_{k}F)-[f_{(z_{q},0)},f_{(0,e_{j})}](e_{k}F)$$
and
$$f_{(x,0)}(E)= f_{(xe_{j},0)}(e_{k}F)-[f_{(x,0)},f_{(0,e_{j})}](e_{k}F).$$

To guarantee that $ \varphi $ preserves Lie multiplication, it remains to prove (3) for $f_{(x,u)}  
$ defined above.

 We apply induction on degree $ n $ of the standard monomial $ E $. If $ n=0 $ (3) is evident. Suppose we have 
 proved (3) for any standard monomial $ E $ of degree less than $ n $ and let $ eE $ be a standard monomial of 
 degree $ n $.
 
We have\
$$f_{([x,y],0)}(eE)= f_{([x,y]e,0)}(E)-\sum_{I*J=E}[f_{([x,y],0)}(I),f_{(0,e)}(J)]=$$
$$f_{([xe,y],0)}(E)+f_{([x,ye],0)}(E) -\sum_{I*J=E}[\sum_{F*H=I}[f_{(x,0)}(F),f_{(y,0)}(H)],f_{(0,e)}(J)] =$$
$$f_{([xe,y],0)}(E) -\sum_{I*J=E}\sum_{F*H=I}[[f_{(x,0)}(F),f_{(0,e)}(J)],f_{(y,0)}(H)] +$$
$$f_{([x,ye],0)}(E) -\sum_{I*J=E}\sum_{F*H=I}[f_{(x,0)}(F),[f_{(y,0)}(H),f_{(0,e)}(J)]] =$$
$$ \sum_{S*H=E}[f_{(xe,0)}(S)-\sum_{F*J=S}[f_{(x,0)}(F),f_{(0,e)}(J)],f_{(y,0)}(H)] +$$
$$ \sum_{F*R=E}[f_{(x,0)}(F),f_{(ye,0)}(R)-\sum_{H*J=R}[f_{(y,0)}(H),f_{(0,e)}(J)]]=$$
$$ \sum_{S*H=E}[(f_{(xe,0)}-[f_{(x,0)},f_{(0,e)}])(S),f_{(y,0)}(H)] +$$ $$ \sum_{F*R=E}[f_{(x,0)}(F),(f_{(ye,0)}-[f_{(y,0)},f_{(0,e)}])(R)]=$$
$$ \sum_{S*H=E}[f_{(x,0)}(eS),f_{(y,0)}(H)] + \sum_{F*R=E}[f_{(x,0)}(F),f_{(y,0)}(eR)]=$$
$$ \sum_{P*Q=eE}[f_{(x,0)}(P),f_{(y,0)}(Q)]= [f_{(x,0)},f_{(y,0)}](eE).$$

The mapping $ \varphi $ is one-to-one. Indeed, $ \varphi((x,u))=\varphi((y,v)) $ implies $ u=v $ and therefore
$f_{(x,u)}=f_{(y,u)}  $. But $f_{(x,u)}=f_{(x,0)}+f_{(0,u)}  $ and $ f_{(y,u)}=f_{(y,0)}+f_{(0,u)} $. Therefore $ f_{(x,0)}= 
f_{(y,0)}$. In particular, $ f_{(x,0)}(1)= f_{(y,0)}(1)$ and $ x=y $.

Thus we have built the mapping which embeds an extension $ N $ of Lie algebra $M  $ by Lie algebra $ L $ into the wreath product $ MWrL $.\\

References\\

1. Kaloujnine L., Sur les $ p $-groupes de Sylow du groups sym\'{e}trique du degr\'{e} $ p^{m} $, C.R. Paris , 1945, 221,     p. 222-224.

2. Kaloujnine L. et Krasner M., Produit complet des groupes de permutations et le probl\`{e}me d'extension des groupes, III, Acta Sci. Math., Szeged, 1951,14, p. 69-82.

3. Jacobson N., Lie Algebras, Interscience Publishers, A division of John Wiley and Sons, New York - London, 1962.

\end{document}